\documentstyle[12pt]{article}

\makeatletter
\renewcommand\thesection{\@Roman\c@section}
\renewcommand\thesubsection{\thesection.\@arabic\c@subsection}
\makeatother

\newcommand{\sect}[1]{\setcounter{equation}{0}\section{#1}}

\newcommand {\beq}{\begin{equation}}
\newcommand {\eeq}{\end{equation}}
\newcommand {\beqa}{\begin{eqnarray}}
\newcommand {\eeqa}{\end{eqnarray}}         %Equation version
\newcommand {\beqs}{\begin{eqnarray*}}
\newcommand {\eeqs}{\end{eqnarray*}}
\newcommand {\bds}{\begin{displaymath}}
\newcommand {\eds}{\end{displaymath}}
\newcommand {\n}{\nonumber\\}

\newcommand{\no}{\noindent}
\newcommand {\bebb}{\begin{thebibliography}}
\newcommand {\eebb}{\end{thebibliography}}      %Reference version
\newcommand {\bbit}{\bibitem}
\def\a{\alpha}

\def\G{\Gamma}

\def\lt{\left}
\def\rt{\right}
\def\rtarr{\rightarrow}

\def\dg{\dagger}

\def\journal#1&#2(#3){\unskip, \sl #1\ \bf #2 \rm(19#3) }
\def\andjournal#1&#2(#3){\sl #1~\bf #2 \rm (19#3) }

\def\npb#1#2#3{Nucl. Phys. {\bf B#1}, #2 (#3)}

\def\plb#1#2#3{Phys. Lett. {\bf B#1}, #2 (#3)}

\def\jmp#1#2#3{J. Math. Phys. {\bf #1}, #2 (#3)}

\begin{document}

%\begin{titlepage}

\begin{flushright}
\end{flushright}

\vskip 1cm

\begin{center}
%\title
{\Large\bf A Unified and Complete Construction of All Finite Dimensional 
Irreducible Representations of $gl(2|2)$} 

\vspace{1cm}

%\author{
{\large Yao-Zhong Zhang and Mark D. Gould}
\vskip.1in
{\em Department of Mathematics, 
University of Queensland, Brisbane, Qld 4072, Australia}

\end{center}

\date{}

%\maketitle

%\vspace{2cm}

\begin{abstract}
Representations of the non-semisimple superalgebra $gl(2|2)$ in the
standard basis are investigated by means of the vector coherent state
method and boson-fermion realization. All finite-dimensional irreducible
typical and atypical representations and lowest weight (indecomposable)
Kac modules of $gl(2|2)$ are constructed explicitly through
the explicit construction of all $gl(2)\oplus gl(2)$ particle states 
(multiplets) in terms of boson and fermion creation operators in the
super-Fock space. This gives a unified and complete treatment of
finite-dimensional representations of $gl(2|2)$ in explicit form,
essential for the construction of primary fields of the corresponding
current superalgebra at arbitrary level.
\end{abstract}

\vspace{1cm}

%%%%%PACS: 11.25Hf; 11.30.Rd; 03.65Fd; 02.20.Hj.

%\vspace{0.5cm}

%\end{titlepage}

\setcounter{section}{0}
\setcounter{equation}{0}
\sect{Introduction}

Recently there is much research interest in superalgebras and their
corresponding non-unitary conformal field theories (CFTs), because of
their applications in high energy and condensed matter physics including
topological field theory \cite{Roz92,Isi94}, logarithmic CFTs (see
e.g. \cite{Flohr} and references therein), disordered systems and the
integer quantum Hall effects
\cite{Efe83,Ber95,Mud96,Maa97,Zir99,Bas00,Gur00,Bha01}. In such
contexts, the vanishing of superdimensions and Virasoro central charges
and the existence of primary fields with negative dimensions are crucial
\cite{Ber95,Mud96}. The most interesting algebras with such properties
are $osp(n|n)$ and $gl(n|n)$.

In most physical applications, one needs the explicit construction of
finite-dimensional representations of a superalgebra. This is
particularly the case in superalgebra CFTs. To construct primary fields
of such CFTs in terms of free fields, one has to construct the
finite-dimensional representations of the superalgebras explicitly.
The explicit construction of the primary fields is essential in
the investigation of disordered systems by the supersymmetric method.

Unlike ordinary bosonic algebras, there are two types of representations
for most superalgebras. They are the so-called typical and atypical
representations. The typical representations are irreducible and are
similar to the usual representations appeared in ordinary bosonic
algebras. The atypical representations have no counterpart in the
bosonic algebra setting. They can be irreducible or not
fully reducible (i.e. reducible or indecomposable). 
This makes the study of representations of superalgebras very difficult.

Representations of $osp(2|2)$ were investigated in \cite{Sch77,Mar80}. A
unified construction of finite-dimensional typical and atypical
representations of
$osp(2|2)$ were given in \cite{Zha03a,Zha03b} by means of the vector
coherent state method. This enabled the explicit construction of
all primary fields of the $osp(2|2)$ CFT \cite{Din03a,Zha03a} in terms of
free fields \cite{Bow96,Ras98}.

In this paper we investigate finite-dimensional
representations of the non-semisimple superalgebra $gl(2|2)$. 
All finite-dimensional irreducible typical and atypical representations 
and lowest weight (indecomposable) Kac modules of $gl(2|2)$
are constructed explicitly through the explicit construction of all
$gl(2)\oplus gl(2)$ particle states (multiplets) in terms of the boson and
fermion creation operators in the super-Fock space. This we believe gives
a unified and complete treatment of all finite-dimensional
irreducible representations of $gl(2|2)$ in explicit form. 

Let us point out that the finite-dimensional representations of $gl(2|2)$
have also been investigated in \cite{Kam89,Pal90} using the GT basis.
Our method is completely different from and in our opinion is simpler
than the method used in these two references. Moreover our results can
be used to construct primary fields of the corresponding $gl(2|2)$ CFTs
at arbitrary level, which is the subject of a separate work. 

This paper is organized as follows. In section 2, we introduce our
notations and derive a free boson-fermion realization of $gl(2|2)$ by
means of the vector coherent state method. In section 3, we describe
the explict construction of independent $gl(2)\oplus gl(2)$ particle states
in the super-Fock space. We derive the actions of odd
simple generators of $gl(2|2)$ on these  multiplets. The 16 independent
multiplets constructed span all finite-dimensional irreducible
typical representations of
$gl(2|2)$. In section 4, we deduce and construct all four types of
finite-dimensional irreducible atypical representations and lowest
weight (indecomposable) Kac modules of $gl(2|2)$.

\setcounter{section}{1}
\setcounter{equation}{0}

\sect{Boson-fermion Realization of $gl(2|2)$}

In this section, we obtain a boson-fermion realization of the
superalgebra $gl(2|2)$ in the standard basis.

This superalgebra is non-semisimple and can be written as 
$gl(2|2)=gl(2|2)^{\rm even}\oplus gl(2|2)^{\rm odd}$, where
\beqa
gl(2|2)^{\rm even}&=&gl(2)\oplus gl(2)\n
&=&\{I\}\oplus\{\{E_{12},E_{21},H_1\}\oplus\{E_{34},E_{43},H_2\},N\},\n
gl(2|2)^{\rm odd}&=&\{E_{13},E_{31},E_{23},E_{32},E_{24},E_{42},
   E_{14},E_{41}\}.
\eeqa
In the standard basis, $E_{12},E_{34},E_{23}$
($E_{21},E_{43},E_{32}$) are simple raising (lowering) generators, 
$E_{13},E_{14},E_{24}$ ($E_{31},E_{41},E_{42}$) are non-simple
raising (lowering) generators and $H_1,H_2,I,N$ are elements of
the Cartan subalgebra. We have  
\beqa
H_1&=&E_{11}-E_{22},~~~~H_2=E_{33}-E_{44},\n
I&=&E_{11}+E_{22}+E_{33}+E_{44},\n
N&=&E_{11}+E_{22}-E_{33}-E_{44}+\beta I\label{cartan-gl22}
\eeqa
with $\beta$ being an arbitrary parameter. That $N$ is not uniquely
determined is a consequence of the fact that $gl(2|2)$ is
non-semisimple. The generators obey the following (anti-)commutation
relations:
\beq
[E_{ij},E_{kl}]=\delta_{jk}E_{il}-(-1)^{([i]+[j])([k]+[l])}\delta_{il}
   E_{kj},
\eeq
where $[E_{ij},E_{kl}]\equiv
E_{ij}E_{kl}-(-1)^{([i]+[j])([k]+[l])}E_{kl}E_{ij}$ is a commutator or
an anticommutator, $[1]=[2]=0,~
[3]=[4]=1$, and $E_{ii},~i=1,2,3,4$ are related to $H_1,H_2,I,N$
via (\ref{cartan-gl22}). The quadratic Casimir of the algebra is given
by $C_2=\sum_{AB}\,(-1)^{[B]}E_{AB}E_{BA}$.

Let $|hw>$ be the highest weight state of highest weight $(J_1,J_2,q,p)$
of $gl(2|2)$ defined by
\beqa
&&H_1|hw>=2J_1|hw>,~~~H_2|hw>=2J_2|hw>,\n
&&I|hw>=2q|hw>,~~~ N|hw>=2p|hw>,\n
&&E_{12}|hw>=E_{34}|hw>=E_{23}|hw>=E_{13}|hw>=E_{14}|hw>=E_{24}|hw>=0.\n
\eeqa
Here $J_1,J_2$ are positive integers and half-integers and $q,p$ are
arbitrary complex numbers.
Define the coherent state \cite{Loh02,Bal88}
$$
e^{E_{21}a_{12}+E_{43}a_{34}+E_{31}\a_{13}+E_{32}\a_{23}+E_{42}\a_{24}
   +E_{41}\a_{14}}|hw>.
$$
Then state vectors are mapped into functions 
\beq
\psi_{J_1,J_2,q,p}=<hw|e^{\a_{13}^\dg E_{13}+\a_{23}^\dg E_{23}
    +\a_{24}^\dg E_{24}+\a_{14}^\dg E_{14}+a_{12}^\dg E_{12}
    +a_{34}^\dg E_{34}}|\psi>|0>,
\eeq 
and operators $A$ are mapped as follows
\beq
A|\psi>\rtarr \G(A)\psi_{J_1,J_2,q,p}=
      <hw|e^{\a_{13}^\dg E_{13}+\a_{23}^\dg E_{23}+\a_{24}^\dg
      E_{24}+\a_{14}^\dg E_{14}+a_{12}^\dg E_{12}
	  +a_{34}^\dg E_{34}}A|\psi>|0>.
\eeq
Here $\a_{ij}^\dg~(\a_{ij})$ are fermion operators with number operators
$N_{\a_{ij}}$ and $a_{ij}^\dg~(a_{ij})$ are boson operators with number
operators $N_{a_{ij}}$. They obey relations
\beqa
&&\{\a_{ij},\a_{kl}^\dg\}=\delta_{ik}\delta_{jl},~~~~\lt(\a_{ij}\rt)^2=
   \lt(\a_{ij}^\dg\rt)^2=0,\n
&&[N_{\a_{ij}},\a_{kl}]=-\delta_{ik}\delta_{jl}\a_{kl},~~~~
  [N_{\a_{ij}},\a_{kl}^\dg]=\delta_{ik}\delta_{jl}\a_{kl}^\dg,\n
&&[a_{ij},a_{kl}^\dg]=\delta_{ik}\delta_{jl},\n
&&[N_{a_{ij}},a_{kl}]=-\delta_{ik}\delta_{jl}a_{kl},~~~~
  [N_{a_{ij}},a_{kl}^\dg]=\delta_{ik}\delta_{jl}a_{kl}^\dg,
\eeqa
and all other (anti-)commutators vanish.
Moreover, $a_{12}|0>=a_{34}|0>=\a_{23}|0>=\a_{13}|0>=\a_{14}|0>=\a_{24}|0>=0$.

Taking $E_{12},E_{34}$ etc in turn and after long algebraic
computations, we find the following representation of simple generators in terms
of the boson and fermion operators:
\beqa
\G(E_{12})&=&a_{12}-\frac{1}{2}\a_{23}^\dg\a_{13}+\frac{1}{2}\lt(\frac{1}{6}
  a_{34}^\dg\a_{23}^\dg-\a_{24}^\dg\rt)\a_{14},\n
\G(E_{34})&=&a_{34}+\frac{1}{2}\a_{23}^\dg\a_{24}+\frac{1}{2}\lt(\frac{1}{6}
  a_{12}^\dg\a_{23}^\dg+\a_{13}^\dg\rt)\a_{14},\n
\G(E_{23})&=&\a_{23}+\frac{1}{2}a_{12}^\dg \a_{13}-\frac{1}{2}a_{34}^\dg
  \lt(\a_{24}+\frac{1}{3}a_{12}^\dg\a_{14}\rt),\n
\G(H_1)&=&2J_1-2N_{a_{12}}+N_{\a_{23}}-N_{\a_{13}}+N_{\a_{24}}-N_{\a_{14}},\n
\G(H_2)&=&2J_2-2N_{a_{34}}+N_{\a_{23}}+N_{\a_{13}}-N_{\a_{24}}-N_{\a_{14}},\n
\G(I)&=&2q,\n
\G(N)&=&2p-2\lt(N_{\a_{23}}+N_{\a_{13}}+N_{\a_{24}}+N_{\a_{14}}\rt),\n
\G(E_{21})&=&a_{12}^\dg\lt[2J_1-N_{a_{12}}+\frac{1}{2}\lt(N_{\a_{23}}-
  N_{\a_{13}}+N_{\a_{24}}-N_{\a_{14}}\rt)\rt]\n
& &-\a_{13}^\dg\a_{23}-\a_{14}^\dg \a_{24}
  -\frac{1}{4}\lt(a_{12}^\dg\rt)^2\a_{23}^\dg\a_{13}\n
& &+\frac{1}{12}a_{12}^\dg a_{34}^\dg\a_{23}^\dg\a_{24}-\frac{1}{4}a_{12}^\dg
   \lt(a_{12}^\dg\a_{24}^\dg+\frac{1}{3}a_{34}^\dg\a_{13}^\dg\rt)\a_{14},\n
\G(E_{43})&=&a_{34}^\dg\lt[2J_2-N_{a_{34}}+\frac{1}{2}\lt(
  N_{\a_{23}}+N_{\a_{13}}-N_{\a_{24}}-N_{\a_{14}}\rt)\rt]\n
& &+\a_{24}^\dg\a_{23}+\a_{14}^\dg\a_{13}+\frac{1}{4}
  \lt(a_{34}^\dg\rt)^2\a_{23}^\dg\a_{24}\n
& &-\frac{1}{12}a_{12}^\dg a_{34}^\dg\a_{23}^\dg\a_{13}
  +\frac{1}{4}\lt(a_{34}^\dg\a_{12}^\dg+\frac{1}{3}a_{12}^\dg\a_{24}^\dg\rt)
  a_{34}^\dg\a_{14},\n
\G(E_{32})&=&\a_{23}^\dg\lt[q-J_1+J_2+\frac{1}{2}(N_{a_{12}}-N_{a_{34}}
  +N_{\a_{13}}-N_{\a_{24}})\rt]\n
& &+\a_{13}^\dg a_{12}+\a_{24}^\dg a_{34}+\frac{1}{6}\a_{23}^\dg
  (a_{12}^\dg\a_{24}^\dg+a_{34}^\dg\a_{13}^\dg)\a_{14},\label{boson-fermion}
\eeqa
and the representation for non-simple generators is easily obtained from that
of simple generators above by means of the commutation relations.
(\ref{boson-fermion}) gives a boson-fermion realization of the non-semisimple
superalgebra $gl(2|2)$ in the standard basis. In this realization,
the Casimir takes a constant value: $C_2=2[(J_1-J_2)(J_1+J_2+1)+q(p-2)]$.

\sect{Typical Representations of $gl(2|2)$}

Representations of $gl(2|2)$ are labelled by $(J_1,J_2,q,p)$ with 
$J_1, J_2$ being positive integers or half-integers and $q,p$ being
arbitrary complex numbers. Consider a particle state in the super-Fock
space,  obtained by acting the
creation operators on the vacuum vector $|0>$.  We call such a state
a level-$x$ state if  $\G(H_1),\G(H_2),\G(I),\G(N)$ have
eigenvalues $2(m_1+x), 2(m_2+x), 2q, 2(p-x)$, respectively. Obviously, a
level-$x$ state is a product of $x$ number of fermion creation operators 
and boson creation operators of the form $\lt(a_{12}^\dg\rt)^{J_1-m_1-y}
\lt(a_{34}^\dg\rt)^{J_2-m_2-\bar{y}}$ acting on $|0>$, 
where $y,\bar{y}$ are certain integers
or half-integers, depending on the values of $x$. 
It is easy to see that there are 16 independent such states obtained
from 16 independent combinations of the
creation operators. This includes one 
level-0 state, four level-1 states, six level-2
states, four level-3 states and one level-4 state. 
Thus each $gl(2|2)$ representation decomposes into
at most 16 representations of the even subalgebra $gl(2)\oplus gl(2)$.
Let us construct representations for $gl(2)\oplus gl(2)$ out of the
above states. First the level-0 and level-4 states are already
representations of $gl(2)\oplus gl(2)$ with highest weights 
$(J_1,J_2,q,p)$ and $(J_1,J_2,q,p-4)$ respectively. We denote these
two multiplets by $|J_1,m_1,J_2,m_2,q;p>$ and $|J_1,m_1,J_2,m_2,q;p-4>$,
respectively. So 
\beqa
&&|J_1,m_1,J_2,m_2,q;p> =\lt(a_{12}^\dg\rt)^{J_1-m_1}\lt(a_{34}^\dg\rt)
  ^{J_2-m_2}|0>,\n
&&~~~~m_1=J_1,J_1-1,\cdots,-J_1,~~m_2=J_2,J_2-1,\cdots, -J_2,\n
&&|J_1,m_1,J_2,m_2,q;p-4>=\a_{23}^\dg\a_{13}^\dg\a_{24}^\dg\a_{14}^\dg
  \lt(a_{12}^\dg\rt)^{J_1-m_1-4}\lt(a_{34}^\dg\rt)^{J_2-m_2-4}|0>,\n
&&~~~~m_1=J_1-4,J_1-5,\cdots,-(J_1+4),~~m_2=J_2-4,J_2-5,\cdots,-(J_2+4).
  \label{level-04b}
\eeqa
Both multiplets have dimension $(2J_1+1)(2J_2+1)$.

It can be shown that other level-$x$ states can be combined into
independent level-$x$ multiplets of $gl(2)\oplus gl(2)$ with certain
highest weights. The procedure is the following. For a given level $x$,
one considers a combination $\Psi_{J_1,m_1,J_2,m_2}$
of all level-$x$ states. The combination
coefficients are in general functions of $J_1,m_1,J_2,m_2$. We
require that $\Psi_{J_1,m_1,J_2,m_2}$ be a representation of
$gl(2)\oplus gl(2)$. In order for the representation to be 
finite-dimensional, the actions
of the $gl(2)\oplus gl(2)$ generators on $\Psi_{J_1,m_1,J_2,m_2}$
must have the following form:
\beqa
&&\G(E_{12})\Psi_{J_1,m_1,J_2,m_2}=(J_1-m_1-z)
    \Psi_{J_1,m_1+1,J_2,m_2},\n
&&\G(E_{21})\Psi_{J_1,m_1,J_2,m_2}=(J_1+m_1+\bar{z})
    \Psi_{J_1,m_1-1,J_2,m_2},\n
&&\G(E_{34})\Psi_{J_1,m_1,J_2,m_2}=(J_2-m_2-u)
    \Psi_{J_1,m_1,J_2,m_2+1},\n
&&\G(E_{43})\Psi_{J_1,m_1,J_2,m_2}=(J_2+m_2+\bar{u})
    \Psi_{J_1,m_1,J_2,m_2-1},
\eeqa
where $z,\bar{z}, u,\bar{u}$ are some integers or half integers to be
determined together with the combination coefficients. These
requirements give rise to difference equations for the combination
coefficients. Solving these difference equations simutaneously for
each level $x$, we determine the combination coefficients and
$z,\bar{z}, u,\bar{u}$. The procedure of solving the difference
equations for each level $x$ is non-trivial and requires long algebraic
manipulations. Here we omit the details and only list the results as
follows.

The four level-1 states can be combined into four independent
multiplets of $gl(2)\oplus gl(2)$ with highest weights $(J_1-\frac{1}{2},
J_2-\frac{1}{2},q,p-1), (J_1+\frac{1}{2},J_2-\frac{1}{2},q,p-1), 
(J_1+\frac{1}{2},J_2+\frac{1}{2},q,p-1)$ and
$(J_1-\frac{1}{2},J_2+\frac{1}{2},q,p-1)$, respectively:
\beqa
&&|J_1-\frac{1}{2},m_1,J_2-\frac{1}{2},m_2,q;p-1>=(\a_{14}^\dg
  +\frac{1}{2}a_{12}^\dg\a_{24}^\dg-\frac{1}{2}\a_{13}^\dg a_{34}^\dg
  -\frac{1}{3}a_{12}^\dg\a_{23}^\dg a_{34}^\dg)\n
&&~~~~\times \lt(a_{12}^\dg\rt)^{J_1-m_1-\frac{3}{2}}\lt(a_{34}^\dg\rt)
  ^{J_2-m_2-\frac{3}{2}}|0>,~~~J_1,J_2\geq\frac{1}{2},\n
&&~~~~m_1=J_1-\frac{3}{2},J_1-\frac{5}{2},\cdots,-(J_1+\frac{1}{2}),~~
  m_2=J_2-\frac{3}{2},J_2-\frac{5}{2},\cdots,-(J_2+\frac{1}{2}),\n
&&|J_1+\frac{1}{2},m_1,J_2-\frac{1}{2},m_2,q;p-1>\n
&&~~=\lt[\frac{1}{2}
  (3J_1+m_1+\frac{5}{2})a_{12}^\dg\a_{24}^\dg-\frac{1}{3}(2J_1+m_1+2)
  a_{12}^\dg\a_{23}^\dg a_{34}^\dg\rt.\n
&&~~~~\lt.-(J_1-m_1-\frac{1}{2})(\a_{14}^\dg-\frac{1}{2}\a_{13}^\dg a_{34}^\dg)
  \rt]\lt(a_{12}^\dg\rt)^{J_1-m_1-\frac{3}{2}}\lt(a_{34}^\dg
  \rt)^{J_2-m_2-\frac{3}{2}}|0>,~~~J_2\geq\frac{1}{2},\n
&&~~~~m_1=J_1-\frac{1}{2},J_1-\frac{3}{2},\cdots,-(J_1+\frac{3}{2}),~~
  m_2=J_2-\frac{3}{2},J_2-\frac{5}{2},\cdots,-(J_2+\frac{1}{2}),\n
&&|J_1+\frac{1}{2},m_1,J_2+\frac{1}{2},m_2,q;p-1>\n
&&~~=\lt[-\frac{1}{4}\lt(
  (3J_1+m_1+\frac{5}{2})(3J_2+m_2+\frac{5}{2})+\frac{1}{3}(J_1-m_1-
  \frac{1}{2})(J_2-m_2-\frac{1}{2})\rt)a_{12}^\dg\a_{23}^\dg a_{34}^\dg\rt.\n
&&~~~~  +\frac{1}{2}(J_1-m_1-\frac{1}{2})(3J_2+m_2+\frac{5}{2})\a_{13}^\dg
  a_{34}^\dg -\frac{1}{2}(3J_1+m_1+\frac{5}{2})(J_2-m_2-\frac{1}{2})a_{12}^\dg
  \a_{24}^\dg\n
&&~~~~\lt.  +(J_1-m_1-\frac{1}{2})(J_2-m_2-\frac{1}{2})\a_{14}^\dg\rt]
  \lt(a_{12}^\dg\rt)^{J_1-m_1-\frac{3}{2}}\lt(a_{34}^\dg
  \rt)^{J_2-m_2-\frac{3}{2}}|0>,\n
&&~~~~m_1=J_1-\frac{1}{2},J_1-\frac{3}{2},\cdots,-(J_1+\frac{3}{2}),~~
  m_2=J_2-\frac{1}{2},J_2-\frac{3}{2},\cdots,-(J_2+\frac{3}{2}),\n
&&|J_1-\frac{1}{2},m_1,J_2+\frac{1}{2},m_2,q;p-1>\n
&&~~=\lt[\frac{1}{2}
  (3J_2+m_2+\frac{5}{2})\a_{13}^\dg a_{34}^\dg+\frac{1}{3}(2J_2+m_2+2)
  a_{12}^\dg\a_{23}^\dg a_{34}^\dg\rt.\n
&&~~~~\lt.+(J_2-m_2-\frac{1}{2})(\a_{14}^\dg-\frac{1}{2}
   a_{12}^\dg\a_{24}^\dg)\rt]
  \lt(a_{12}^\dg\rt)^{J_1-m_1-\frac{3}{2}}\lt(a_{34}^\dg
  \rt)^{J_2-m_2-\frac{3}{2}}|0>,~~~J_1\geq\frac{1}{2},\n
&&~~~~m_1=J_1-\frac{3}{2},J_1-\frac{5}{2},\cdots,-(J_1+\frac{1}{2}),~~
  m_2=J_2-\frac{1}{2},J_2-\frac{3}{2},\cdots,-(J_2+\frac{3}{2}).
  \label{level-1b}
\eeqa
The dimensions for these multiplets are $(2J_1)(2J_2),
(2J_1+2)(2J_2), (2J_1+2)(2J_2+2)$ and $(2J_1) (2J_2+2)$, respectively.

The six level-2 states can be combined into 6 independent multiplets of
$gl(2)\oplus gl(2)$ with highest weights $(J_1,J_2-1,q,p-2),
(J_1-1,J_2,q,p-2), (J_1+1,J_2,q,p-2), (J_1,J_2+1,q,p-2),
(J_1,J_2,q,p-2)$ and $(J_1,J_2,q,p-2)$, respectively:
\beqa
&&|J_1,m_1,J_2-1,m_2,q;p-2>=\a_{24}^\dg\a_{14}^\dg
  \lt(a_{12}^\dg\rt)^{J_1-m_1-2}\lt(a_{34}^\dg\rt)^{J_2-m_2-3}|0>\n
&&~~~~+\frac{1}{2}\lt[-\a_{23}^\dg\a_{14}^\dg
  +\frac{1}{6}\a_{23}^\dg\a_{24}^\dg a_{12}^\dg+\a_{13}^\dg\a_{24}^\dg
  +\frac{1}{2}\a_{23}^\dg\a_{13}^\dg a_{34}^\dg\rt]
  \lt(a_{12}^\dg\rt)^{J_1-m_1-2}\lt(a_{34}^\dg\rt)^{J_2-m_2-2}|0>,\n
&&~~~~J_2\geq 1,~~~m_1=J_1-2,J_1-3,\cdots, -(J_1+2),~~
  m_2=J_2-3,J_2-4,\cdots,-(J_2+1),\n
&&|J_1-1,m_1,J_2,m_2,q;p-2>=\a_{13}^\dg\a_{14}^\dg
  \lt(a_{12}^\dg\rt)^{J_1-m_1-3}\lt(a_{34}^\dg\rt)^{J_2-m_2-2}|0>\n
&&~~~~+\frac{1}{2}\lt[\a_{23}^\dg\a_{14}^\dg
  +\frac{1}{6}\a_{23}^\dg\a_{13}^\dg a_{34}^\dg+\a_{13}^\dg\a_{24}^\dg
  +\frac{1}{2}\a_{23}^\dg\a_{24}^\dg a_{12}^\dg\rt]
  \lt(a_{12}^\dg\rt)^{J_1-m_1-2}\lt(a_{34}^\dg\rt)^{J_2-m_2-2}|0>,\n
&&~~~~J_1\geq 1,~~~m_1=J_1-3,J_1-4,\cdots, -(J_1+1),~~  
  m_2=J_2-2,J_2-3,\cdots,-(J_2+2),\n
&&|J_1+1,m_1,J_2,m_2,q;p-2>\n
&&~~=\lt[\frac{1}{2}[J_1-m_1-1+(3J_1+m_1+3)(3J_1+m_1+5)]
  \a_{23}^\dg\a_{24}^\dg \lt(a_{12}^\dg\rt)^2\rt.\n
&&~~~~  +(J_1-m_1-1)(J_1-m_1-2)(\a_{13}^\dg\a_{14}^\dg+
  \frac{1}{12}a_{12}^\dg\a_{23}^\dg\a_{13}^\dg a_{34}^\dg)\n
&&~~~~\lt.-\frac{1}{2}(J_1-m_1-1)(3J_1+m_1+4)a_{12}^\dg(\a_{13}^\dg\a_{24}^\dg
  +\a_{23}^\dg\a_{14}^\dg)\rt] \lt(a_{12}^\dg\rt)^{J_1-m_1-3}
  \lt(a_{34}^\dg\rt)^{J_2-m_2-2}|0>,\n
&&~~~~m_1=J_1-1,J_1-2,\cdots, -(J_1+3),~~m_2=J_2-2,J_2-3,\cdots,-(J_2+2),\n
&&|J_1,m_1,J_2+1,m_2,q;p-2>\n
&&~~=\lt[\frac{1}{4}[J_2-m_2-1+(3J_2+m_2+3)(3J_2+m_2+5)]
  \a_{23}^\dg\a_{13}^\dg\lt(a_{34}^\dg\rt)^2\rt.\n
&&~~~~+\frac{1}{2}(J_2-m_2-1)(3J_2+m_2+4)(\a_{23}^\dg\a_{14}^\dg
  -\a_{13}^\dg\a_{24}^\dg)a_{34}^\dg\n
&&~~~~\lt.+(J_2-m_2-1)(J_2-m_2-2)(\a_{24}^\dg\a_{14}^\dg
  +\frac{1}{12}\a_{23}^\dg\a_{24}^\dg a_{12}^\dg a_{34}^\dg)\rt]\n
&&~~~~\times\lt(a_{12}^\dg\rt)^{J_1-m_1-2}\lt(a_{34}^\dg\rt)^{J_2-m_2-3}|0>,\n
&&~~~~m_1=J_1-2,J_1-3,\cdots, -(J_1+2),~~m_2=J_2-1,J_2-2,\cdots,-(J_2+3),\n
&&|J_1,m_1,J_2,m_2,q;p-2>_{\bf I}=(J_2-m_2-2)\a_{24}^\dg\a_{14}^\dg
  \lt(a_{12}^\dg\rt)^{J_1-m_1-2}\lt(a_{34}^\dg\rt)^{J_2-m_2-3}|0>\n
&&~~~~+\lt[\frac{1}{2}(J_2+m_2+2)(\a_{23}^\dg\a_{14}^\dg-\a_{13}^\dg\a_{24}^\dg)
  +\frac{1}{12}(J_2-m_2-2)\a_{23}^\dg\a_{24}^\dg a_{12}^\dg\rt.\n
&&~~~~\lt.-\frac{1}{4}(3J_2+m_2+2)\a_{23}^\dg\a_{13}^\dg a_{34}^\dg\rt]
  \lt(a_{12}^\dg\rt)^{J_1-m_1-2}\lt(a_{34}^\dg\rt)^{J_2-m_2-2}|0>,\n
&&~~~~m_1=J_1-2,J_1-3,\cdots,-(J_1+2),~~m_2=J_2-1,J_2-2,\cdots,-(J_2+2),\n
&&|J_1,m_1,J_2,m_2,q;p-2>_{\bf II}=(J_1-m_1-2)\a_{13}^\dg\a_{14}^\dg
  \lt(a_{12}^\dg\rt)^{J_1-m_1-3}\lt(a_{34}^\dg\rt)^{J_2-m_2-2}|0>\n
&&~~~~+\lt[-\frac{1}{2}(J_1+m_1+2)(\a_{13}^\dg\a_{24}^\dg
  +\a_{23}^\dg\a_{14}^\dg)
  +\frac{1}{12}(J_1-m_1-2)\a_{23}^\dg\a_{13}^\dg a_{34}^\dg\rt.\n
&&~~~~\lt.-\frac{1}{4}(3J_1+m_1+2)\a_{23}^\dg\a_{24}^\dg a_{12}^\dg\rt]
  \lt(a_{12}^\dg\rt)^{J_1-m_1-2}\lt(a_{34}^\dg\rt)^{J_2-m_2-2}|0>,\n
&&~~~~m_1=J_1-2,J_1-3,\cdots,-(J_1+2),~~m_2=J_2-1,J_2-2,\cdots,-(J_2+2).
  \label{level-2b}
\eeqa
Notice that the last two multiplets, which have been denoted above by
$|J_1,m_1,J_2,m_2,q; p-2>_{\bf I}$ and $|J_1,m_1,J_2,m_2,q; p-2>_{\bf II}$,
respectively, have the same highest weight $(J_1,J_2,q,p-2)$. 
This means that multiplicity will in general appear in the
$gl(2|2)\downarrow gl(2)\oplus gl(2)$ branching rule.
It is easy to see from the above expressions
that $|J_1,m_1,J_2,m_2,q;p-2>_{\bf I}\equiv 0$ when
$J_2=0$ and $|J_1,m_1,J_2,m_2,q;p-2>_{\bf II}\equiv 0$ when $J_1=0$.

The dimensions for the first four multiplets are $(2J_1+1)(2J_2-1),
(2J_1-1)(2J_2+1), (2J_1+3)(2J_2+1)$ and $(2J_1+1)(2J_2+3)$, respectively. 
The dimension for $|J_1,m_1,J_2,m_2,q;p-2>_{\bf I}$ is
$(2J_1+1)(2J_2+1)$ if $J_2\neq 0$ and zero if $J_2=0$. Similarly,
the dimension for $|J_1,m_1,J_2,m_2,q;p-2>_{\bf II}$ is
$(2J_1+1)(2J_2+1)$ if $J_1\neq 0$ and zero if $J_1=0$.

Finally, the four level-3 states are combined into four independent
multiplets of $gl(2)\oplus gl(2)$ with highest weights $(J_1-\frac{1}{2},
J_2-\frac{1}{2},q,p-3), (J_1+\frac{1}{2},J_2-\frac{1}{2},q,p-3),
(J_1-\frac{1}{2},J_2+\frac{1}{2},q,p-3)$ and
$(J_1+\frac{1}{2},J_2+\frac{1}{2},q,p-3)$, respectively:
\beqa
&&|J_1-\frac{1}{2},m_1,J_2-\frac{1}{2},m_2,q;p-3>
  =\lt[(\a_{13}^\dg+\frac{1}{2}a_{12}^\dg\a_{23}^\dg)
  \a_{24}^\dg\a_{14}^\dg\rt.\n
&&~~~~\lt.+\frac{1}{2}\a_{23}^\dg\a_{13}^\dg(\a_{14}^\dg a_{34}^\dg
  +\frac{1}{3}a_{12}^\dg\a_{24}^\dg a_{34}^\dg)\rt]
  \lt(a_{12}^\dg\rt)^{J_1-m_1-\frac{7}{2}}\lt(a_{34}^\dg\rt)^{J_2-m_2-
  \frac{7}{2}}|0>,\n
&&~~~~J_1,J_2\geq \frac{1}{2},~~~
  m_1=J_1-\frac{7}{2},\cdots,-(J_1+\frac{5}{2}),~~
  m_2=J_2-\frac{7}{2},\cdots,-(J_2+\frac{5}{2}),\n
&&|J_1+\frac{1}{2},m_1,J_2-\frac{1}{2},m_2,q;p-3>
  =\lt[-\frac{1}{2}(3J_1+m_1+\frac{9}{2})\a_{23}^\dg\a_{24}^\dg\a_{14}^\dg
   a_{12}^\dg\rt.\n
&&~~~~+(J_1-m_1-\frac{5}{2})(-\a_{24}^\dg+\frac{1}{2}\a_{23}^\dg 
  a_{34}^\dg)\a_{13}^\dg\a_{14}^\dg\n
&&~~~~\lt.-\frac{1}{6}(5J_1+m_1+\frac{11}{2})\a_{23}^\dg\a_{13}^\dg
  \a_{24}^\dg a_{12}^\dg a_{34}^\dg\rt]
  \lt(a_{12}^\dg\rt)^{J_1-m_1-\frac{7}{2}}\lt(a_{34}^\dg\rt)^{J_2-m_2-
  \frac{7}{2}}|0>,\n
&&~~~~J_2\geq \frac{1}{2},~~~  
  m_1=J_1-\frac{5}{2},\cdots,-(J_1+\frac{7}{2}),~~
    m_2=J_2-\frac{7}{2},\cdots,-(J_2+\frac{5}{2}),\n
&&|J_1-\frac{1}{2},m_1,J_2+\frac{1}{2},m_2,q;p-3>
  =\lt[-\frac{1}{2}(3J_2+m_2+\frac{9}{2})\a_{23}^\dg\a_{13}^\dg\a_{14}^\dg
   a_{34}^\dg\rt.\n
&&~~~~+(J_2-m_2-\frac{5}{2})(\a_{13}^\dg+\frac{1}{2}\a_{23}^\dg
   a_{12}^\dg)\a_{24}^\dg\a_{14}^\dg\n
&&~~~~\lt.-\frac{1}{6}(5J_2+m_2+\frac{11}{2})\a_{23}^\dg\a_{13}^\dg  
  \a_{24}^\dg a_{12}^\dg a_{34}^\dg\rt]
  \lt(a_{12}^\dg\rt)^{J_1-m_1-\frac{7}{2}}\lt(a_{34}^\dg\rt)^{J_2-m_2-
  \frac{7}{2}}|0>,\n
&&~~~~J_1\geq \frac{1}{2},~~~  
  m_1=J_1-\frac{7}{2},\cdots,-(J_1+\frac{5}{2}),~~
  m_2=J_2-\frac{5}{2},\cdots,-(J_2+\frac{7}{2}),\n
&&|J_1+\frac{1}{2},m_1,J_2+\frac{1}{2},m_2,q;p-3>
  =\lt[\frac{1}{4}\lt((3J_1+m_1+\frac{9}{2})(3J_2+m_2+\frac{9}{2})\rt.\rt.\n
&&~~~~\lt.  -\frac{1}{3}(J_1-m_1-\frac{5}{2})(J_2-m_2-\frac{5}{2})\rt)
  a_{12}^\dg\a_{23}^\dg\a_{13}^\dg\a_{24}^\dg a_{34}^\dg\n
&&~~~~-\frac{1}{2}(J_1-m_1-\frac{5}{2})(3J_2+m_2+\frac{9}{2})
  \a_{23}^\dg\a_{13}^\dg\a_{14}^\dg a_{34}^\dg\n
&&~~~~-\frac{1}{2}(3J_1+m_1+\frac{9}{2})(J_2-m_2-\frac{5}{2})
  a_{12}^\dg\a_{23}^\dg\a_{24}^\dg\a_{14}^\dg\n
&&~~~~\lt.+(J_1-m_1-\frac{5}{2})(J_2-m_2-\frac{5}{2})\a_{13}^\dg
  \a_{24}^\dg\a_{14}^\dg\rt]
  \lt(a_{12}^\dg\rt)^{J_1-m_1-\frac{7}{2}}\lt(a_{34}^\dg\rt)^{J_2-m_2-
  \frac{7}{2}}|0>,\n
&&~~~~m_1=J_1-\frac{5}{2},\cdots,-(J_1+\frac{7}{2}),~~
      m_2=J_2-\frac{5}{2},\cdots,-(J_2+\frac{7}{2}).
      \label{level-3b}
\eeqa
The dimensions for these multiplets are $(2J_1)(2J_2),
(2J_1+2)(2J_2), (2J_1)(2J_2+2)$ and $(2J_1+2)(2J_2+2)$, respectively.

The actions of the odd generators of $gl(2|2)$ on the $gl(2)\oplus gl(2)$
multiplets (\ref{level-04b}) and (\ref{level-1b})-(\ref{level-3b})
can be computed by means of the free boson-fermion
realization of the generators. 
In the following we list the actions of the odd simple
generators. The actions of odd non-simple generators can be easily obtained
using the commutation relations.

First for the level-0 multiplet, we have the actions of the odd simple
generators
\beqa
&&\G(E_{23})|J_1,m_1,J_2,m_2,q;p>=0,\n
&&\G(E_{32})|J_1,m_1,J_2,m_2,q;p>=\frac{1}{(2J_1+1)(2J_2+1)}\n
&&~~~~\times\lt[-(q+J_1-J_2)(J_1-m_1)(J_2-m_2)
  |J_1-\frac{1}{2},m_1-\frac{1}{2},J_2-\frac{1}{2},
  m_2-\frac{1}{2},q;p-1> \rt.\n
&&~~~~-(q-J_1-J_2-1)(J_2-m_2)|J_1+\frac{1}{2},m_1-\frac{1}{2},J_2-\frac{1}{2},
  m_2-\frac{1}{2},q;p-1>\n
&&~~~~-(q-J_1+J_2)|J_1+\frac{1}{2},m_1-\frac{1}{2},J_2+\frac{1}{2},
  m_2-\frac{1}{2},q;p-1>\n
&&~~~~\lt.+(q+J_1+J_2+1)(J_1-m_1)
  |J_1-\frac{1}{2},m_1-\frac{1}{2},J_2+\frac{1}{2},
  m_2-\frac{1}{2},q;p-1>\rt].\label{action-level0}
\eeqa
{}From (\ref{action-level0}) we see that when $q=J_1-J_2$ (resp. $-J_1+J_2$)
the third (resp. first) term vanishes and if $q=J_1+J_2+1$ (resp. 
$-J_1-J_2-1$) then
the second (resp. fourth) term disappears. This indicates that when
$q=\pm(J_1-J_2), \pm(J_1+J_2+1)$ atypical representations arise (see
next section for details).

For the four level-1 multiplets, we obtain the the following actions of the odd
simple generators, after long algebraic manipulations,
\beqa
&&\G(E_{23})|J_1-\frac{1}{2},m_1,J_2-\frac{1}{2},m_2,q,p-1>=
  -|J_1,m_1+\frac{1}{2},J_2,m_2+\frac{1}{2},q;p>,\n
&&\G(E_{23})|J_1+\frac{1}{2},m_1,J_2-\frac{1}{2},m_2,q,p-1>\n
&&~~= -(J_1+m_1+\frac{3}{2})|J_1,m_1+\frac{1}{2},J_2,m_2+\frac{1}{2},q;p>,\n
&&\G(E_{23})|J_1+\frac{1}{2},m_1,J_2+\frac{1}{2},m_2,q,p-1>\n
&&~~= -(J_1+m_1+\frac{3}{2})(J_2+m_2+\frac{3}{2})
  |J_1,m_1+\frac{1}{2},J_2,m_2+\frac{1}{2},q;p>,\n
&&\G(E_{23})|J_1-\frac{1}{2},m_1,J_2+\frac{1}{2},m_2,q,p-1>\n
&&~~= -(J_2+m_2+\frac{3}{2})|J_1,m_1+\frac{1}{2},J_2,m_2+\frac{1}{2},q;p>,\n
&&\G(E_{32})|J_1-\frac{1}{2},m_1,J_2-\frac{1}{2},m_2,q,p-1>\n
&&~~= -\frac{J_2-m_2-\frac{3}{2}}{2J_2}(q-J_1-J_2-1)
  |J_1,m_1-\frac{1}{2},J_2-1,m_2-\frac{1}{2},q;p-2>\n
&&~~~~+\frac{q-J_1+J_2-1}{2J_2}|J_1,m_1-\frac{1}{2},J_2,m_2-\frac{1}{2},
  q;p-2>_{\bf I}\n
&&~~~~-\frac{q-J_1+J_2+1}{2J_2}|J_1,m_1-\frac{1}{2},J_2,m_2-\frac{1}{2},
  q;p-2>_{\bf II}\n
&&~~~~+\frac{J_1-m_1-\frac{3}{2}}{2J_1}(q+J_1+J_2+1)
  |J_1-1,m_1-\frac{1}{2},J_2,m_2-\frac{1}{2},q;p-2>,\n
&&\G(E_{32})|J_1+\frac{1}{2},m_1,J_2-\frac{1}{2},m_2,q,p-1>\n
&&~~=  \frac{(J_1-m_1-\frac{1}{2})(J_2-m_2-\frac{3}{2})}{2J_2}(q+J_1-J_2)
  |J_1,m_1-\frac{1}{2},J_2-1,m_2-\frac{1}{2},q;p-2>\n
&&~~~~-\frac{J_1-m_1-\frac{1}{2}}{2J_2}(q+J_1+J_2)
  |J_1,m_1-\frac{1}{2},J_2,m_2-\frac{1}{2}, q;p-2>_{\bf I}\n
&&~~~~-\frac{J_1-m_1-\frac{1}{2}}{2(J_1+1)}(q+J_1+J_2+2)
  |J_1,m_1-\frac{1}{2},J_2,m_2-\frac{1}{2},q;p-2>_{\bf II}\n
&&~~~~+\frac{q-J_1+J_2}{2(J_1+1)}  
  |J_1+1,m_1-\frac{1}{2},J_2,m_2-\frac{1}{2},q;p-2>,\n
&&\G(E_{32})|J_1+\frac{1}{2},m_1,J_2+\frac{1}{2},m_2,q,p-1>\n
&&~~= -\frac{(J_1-m_1-\frac{1}{2})(J_2-m_2-\frac{1}{2})}{2(J_2+1)}
  (q+J_1-J_2-1) |J_1,m_1-\frac{1}{2},J_2,m_2-\frac{1}{2},q;p-2>_{\bf I}\n
&&~~~~+\frac{(J_1-m_1-\frac{1}{2})(J_2-m_2-\frac{1}{2})}{2(J_1+1)}(q+J_1-J_2+1)
      |J_1,m_1-\frac{1}{2},J_2,m_2-\frac{1}{2}, q;p-2>_{\bf II}\n
&&~~~~-\frac{J_2-m_2-\frac{1}{2}}{2(J_1+1)}(q-J_1-J_2-1)
  |J_1+1,m_1-\frac{1}{2},J_2,m_2-\frac{1}{2},q;p-2>\n
&&~~~~+\frac{J_1-m_1-\frac{1}{2}}{2(J_2+1)}(q+J_1+J_2+1) 
  |J_1,m_1-\frac{1}{2},J_2+1,m_2-\frac{1}{2},q;p-2>,\n
&&\G(E_{32})|J_1-\frac{1}{2},m_1,J_2+\frac{1}{2},m_2,q,p-1>\n
&&~~=  -\frac{J_2-m_2-\frac{1}{2}}{2(J_2+1)}(q-J_1-J_2-2)
  |J_1,m_1-\frac{1}{2},J_2,m_2-\frac{1}{2},q;p-2>_{\bf I}\n
&&~~~~-\frac{J_2-m_2-\frac{1}{2}}{2J_1}(q-J_1-J_2)
  |J_1,m_1-\frac{1}{2},J_2,m_2-\frac{1}{2}, q;p-2>_{\bf II}\n
&&~~~~+\frac{(J_1-m_1-\frac{3}{2})(J_2-m_2-\frac{1}{2})}{2J_1}(q+J_1-J_2)
  |J_1-1,m_1-\frac{1}{2},J_2,m_2-\frac{1}{2},q;p-2>\n
&&~~~~+\frac{q-J_1+J_2}{2(J_2+1)} 
  |J_1,m_1-\frac{1}{2},J_2+1,m_2-\frac{1}{2},q;p-2>.\label{action-level1}
\eeqa

Similar to the level-1 case, we find after long algebraic computations
that the actions of the odd simple generators on the six
level-2 multiplets are given by
\beqa
&&\G(E_{23})|J_1,m_1,J_2-1,m_2,q;p-2>\n
&&~~=-\frac{J_1+m_1+2}{2J_1+1}|J_1-\frac{1}{2},m_1+\frac{1}{2},J_2-\frac{1}{2},
  m_2+\frac{1}{2},q,p-1>\n
&&~~~~+\frac{1}{2J_1+1}|J_1+\frac{1}{2},m_1+\frac{1}{2},J_2-\frac{1}{2},
  m_2+\frac{1}{2},q,p-1>,\n
&&\G(E_{23})|J_1-1,m_1,J_2,m_2,q;p-2>\n
&&~~=\frac{J_2+m_2+2}{2J_2+1}|J_1-\frac{1}{2},m_1+\frac{1}{2},J_2-\frac{1}{2},
  m_2+\frac{1}{2},q,p-1>\n
&&~~~~+\frac{1}{2J_1+1}|J_1-\frac{1}{2},m_1+\frac{1}{2},J_2+\frac{1}{2},  
  m_2+\frac{1}{2},q,p-1>,\n
&&\G(E_{23})|J_1+1,m_1,J_2,m_2,q;p-2>\n
&&~~=\frac{(J_1+m_1+3)(J_2+m_2+2)}{2J_2+1}
  |J_1+\frac{1}{2},m_1+\frac{1}{2},J_2-\frac{1}{2}, m_2+\frac{1}{2},q,p-1>\n
&&~~~~-\frac{J_1+m_1+3}{2J_2+1}|J_1+\frac{1}{2},m_1+\frac{1}{2},
  J_2+\frac{1}{2}, m_2+\frac{1}{2},q,p-1>,\n
&&\G(E_{23})|J_1,m_1,J_2+1,m_2,q;p-2>\n  
&&~~=\frac{J_2+m_2+3}{2J_1+1}
  |J_1+\frac{1}{2},m_1+\frac{1}{2},J_2+\frac{1}{2}, m_2+\frac{1}{2},q,p-1>\n
&&~~~~+\frac{(J_1+m_1+2)(J_2+m_2+3)}{2J_1+1}|J_1-\frac{1}{2},m_1+\frac{1}{2},
    J_2+\frac{1}{2}, m_2+\frac{1}{2},q,p-1>,\n
&&\G(E_{23})|J_1,m_1,J_2,m_2,q;p-2>_{\bf I}=\frac{1}{(2J_1+1)(2J_2+1)}\times\n
&&~~~~\lt[(J_2+1)(J_1+m_1+2)(J_2+m_2+2)|J_1-\frac{1}{2},m_1+\frac{1}{2},
  J_2-\frac{1}{2},m_2+\frac{1}{2},q,p-1>\rt.\n
&&~~~~-(J_2+1)(J_2+m_2+2)|J_1+\frac{1}{2},m_1+\frac{1}{2},
  J_2-\frac{1}{2},m_2+\frac{1}{2},q,p-1>\n
&&~~~~-J_2|J_1+\frac{1}{2},m_1+\frac{1}{2},J_2+\frac{1}{2},
  m_2+\frac{1}{2},q,p-1>\n
&&~~~~\lt.-J_2(J_1+m_1+2)|J_1-\frac{1}{2},m_1+\frac{1}{2},J_2+\frac{1}{2},
  m_2+\frac{1}{2},q,p-1>\rt],\n
&&\G(E_{23})|J_1,m_1,J_2,m_2,q;p-2>_{\bf II}=\frac{1}{(2J_1+1)(2J_2+1)}\times\n
&&~~~~\lt[-(J_1+1)(J_1+m_1+2)(J_2+m_2+2)|J_1-\frac{1}{2},m_1+\frac{1}{2},
 J_2-\frac{1}{2},m_2+\frac{1}{2},q,p-1>\rt.\n
&&~~~~-J_1(J_2+m_2+2)|J_1+\frac{1}{2},m_1+\frac{1}{2},  
   J_2-\frac{1}{2},m_2+\frac{1}{2},q,p-1>\n
&&~~~~+J_1|J_1+\frac{1}{2},m_1+\frac{1}{2},J_2+\frac{1}{2},
   m_2+\frac{1}{2},q,p-1>\n
&&~~~~\lt.-(J_1+1)(J_1+m_1+2)|J_1-\frac{1}{2},m_1+\frac{1}{2},J_2+\frac{1}{2},
   m_2+\frac{1}{2},q,p-1>\rt],\n
&&\G(E_{32})|J_1,m_1,J_2-1,m_2,q;p-2>\n
&&~~=\frac{J_1-m_1-2}{2J_1+1}(q+J_1+J_2+1)
  |J_1-\frac{1}{2},m_1-\frac{1}{2},J_2-\frac{1}{2},m_2-\frac{1}{2},q;p-3>\n
&&~~~~-\frac{q-J_1+J_2}{2J_1+1}
  |J_1+\frac{1}{2},m_1-\frac{1}{2},J_2-\frac{1}{2},m_2-\frac{1}{2},q;p-3>,\n
&&\G(E_{32})|J_1-1,m_1,J_2,m_2,q;p-2>\n
&&~~=\frac{J_2-m_2-2}{2J_2+1}(q-J_1-J_2-1)
  |J_1-\frac{1}{2},m_1-\frac{1}{2},J_2-\frac{1}{2},m_2-\frac{1}{2},q;p-3>\n
&&~~~~-\frac{q-J_1+J_2}{2J_2+1}
  |J_1-\frac{1}{2},m_1-\frac{1}{2},J_2+\frac{1}{2},m_2-\frac{1}{2},q;p-3>,\n
&&\G(E_{32})|J_1+1,m_1,J_2,m_2,q;p-2> 
  =\frac{(J_1-m_1-1)(J_2-m_2-2)}{2J_2+1}\n
&&~~~~\times(q+J_1-J_2)
  |J_1+\frac{1}{2},m_1-\frac{1}{2},J_2-\frac{1}{2},m_2-\frac{1}{2},q;p-3>\n
&&~~~~-\frac{J_1-m_1-1}{2J_2+1}(q+J_1+J_2+1)
  |J_1+\frac{1}{2},m_1-\frac{1}{2},J_2+\frac{1}{2},m_2-\frac{1}{2},q;p-3>,\n
&&\G(E_{32})|J_1,m_1,J_2+1,m_2,q;p-2>
   =\frac{(J_1-m_1-2)(J_2-m_2-1)}{2J_1+1}\n
&&~~~~\times(q+J_1-J_2)
  |J_1-\frac{1}{2},m_1-\frac{1}{2},J_2+\frac{1}{2},m_2-\frac{1}{2},q;p-3>\n
&&~~~~-\frac{J_2-m_2-1}{2J_1+1}(q-J_1-J_2-1)
  |J_1+\frac{1}{2},m_1-\frac{1}{2},J_2+\frac{1}{2},m_2-\frac{1}{2},q;p-3>,\n
&&\G(E_{32})|J_1,m_1,J_2,m_2,q;p-2>_{\bf I}=\frac{1}{(2J_1+1)(2J_2+1)}
   \lt[(J_2+1)(J_1-m_1-2)\rt.\n
&&~~~~ \times(J_2-m_2-2)(q+J_1-J_2+1)
  |J_1-\frac{1}{2},m_1-\frac{1}{2},J_2-\frac{1}{2},m_2-\frac{1}{2},q;p-3>\n
&&~~~~-(J_2+1)(J_2-m_2-2)(q-J_1-J_2)
  |J_1+\frac{1}{2},m_1-\frac{1}{2},J_2-\frac{1}{2},m_2-\frac{1}{2},q;p-3>\n
&&~~~~+J_2(J_1-m_1-2)(q+J_1+J_2+2)
  |J_1-\frac{1}{2},m_1-\frac{1}{2},J_2+\frac{1}{2},m_2-\frac{1}{2},q;p-3>\n
&&~~~~\lt.-J_2(q-J_1+J_2+1)
 |J_1+\frac{1}{2},m_1-\frac{1}{2},J_2+\frac{1}{2},m_2-\frac{1}{2},q;p-3>\rt],\n
&&\G(E_{32})|J_1,m_1,J_2,m_2,q;p-2>_{\bf II}=\lt.\frac{1}{(2J_1+1)(2J_2+1)}
  \rt[(J_1+1)(J_1-m_1-2)\n
&&~~~~\times(J_2-m_2-2)(q+J_1-J_2-1)
  |J_1-\frac{1}{2},m_1-\frac{1}{2},J_2-\frac{1}{2},m_2-\frac{1}{2},q;p-3>\n
&&~~~~+J_1(J_2-m_2-2)(q-J_1-J_2-2)
  |J_1+\frac{1}{2},m_1-\frac{1}{2},J_2-\frac{1}{2},m_2-\frac{1}{2},q;p-3>\n
&&~~~~-(J_1+1)(J_1-m_1-2)(q+J_1+J_2)
  |J_1-\frac{1}{2},m_1-\frac{1}{2},J_2+\frac{1}{2},m_2-\frac{1}{2},q;p-3>\n
&&~~~~\lt.-J_1(q-J_1+J_2-1)
 |J_1+\frac{1}{2},m_1-\frac{1}{2},J_2+\frac{1}{2},m_2-\frac{1}{2},q;p-3>\rt].
 \label{action-level2}
\eeqa

The actions of the odd simple generators on the four level-3 multiplets
can be obtained in a similar way. We list the results as follows:
\beqa
&&\G(E_{23})|J_1-\frac{1}{2},m_1,J_2-\frac{1}{2},m_2,q,p-3>\n
&&~~= \frac{J_2+m_2+\frac{5}{2}}{2J_2}
   |J_1,m_1+\frac{1}{2},J_2-1,m_2+\frac{1}{2},q;p-2>\n
&&~~~~+\frac{1}{2J_2}|J_1,m_1+\frac{1}{2},J_2,m_2+\frac{1}{2},q;p-2>_{\bf I}
   +\frac{1}{2J_1}|J_1,m_1+\frac{1}{2},J_2,m_2+\frac{1}{2},q;p-2>_{\bf II}\n
&&~~~~+\frac{J_1+m_1+\frac{5}{2}}{2J_1}
   |J_1-1,m_1+\frac{1}{2},J_2,m_2+\frac{1}{2},q;p-2>,\n
&&\G(E_{23})|J_1+\frac{1}{2},m_1,J_2-\frac{1}{2},m_2,q,p-3>\n
&&~~=-\frac{(J_1+m_1+\frac{7}{2})(J_2+m_2+\frac{5}{2})}{2J_2}
  |J_1,m_1+\frac{1}{2},J_2-1,m_2+\frac{1}{2},q;p-2>\n
&&~~~~+(J_1+m_1+\frac{7}{2})\lt[-\frac{1}{2J_2}|J_1,m_1+\frac{1}{2},J_2,
  m_2+\frac{1}{2},q;p-2>_{\bf I}\rt.\n
&&~~~~\lt.+\frac{1}{2(J_1+1)}
  |J_1,m_1+\frac{1}{2},J_2,m_2+\frac{1}{2},q;p-2>_{\bf II}\rt]\n
&&~~~~+\frac{1}{2(J_1+1)}|J_1+1,m_1+\frac{1}{2},J_2,m_2+\frac{1}{2},q;p-2>,\n
&&\G(E_{23})|J_1-\frac{1}{2},m_1,J_2+\frac{1}{2},m_2,q,p-3>\n
&&~~=(J_2+m_2+\frac{7}{2})\lt[\frac{1}{2(J_2+1)}|J_1,m_1+\frac{1}{2},J_2,
  m_2+\frac{1}{2},q;p-2>_{\bf I}\rt.\n
&&~~~~\lt.-\frac{1}{2J_1}
  |J_1,m_1+\frac{1}{2},J_2,m_2+\frac{1}{2},q;p-2>_{\bf II}\rt]\n
&&~~~~-\frac{(J_1+m_1+\frac{5}{2})(J_2+m_2+\frac{7}{2})}{2J_1}
  |J_1-1,m_1+\frac{1}{2},J_2,m_2+\frac{1}{2},q;p-2>\n
&&~~~~+\frac{1}{2(J_2+1)}|J_1,m_1+\frac{1}{2},J_2+1,m_2+\frac{1}{2},q;p-2>,\n
&&\G(E_{23})|J_1+\frac{1}{2},m_1,J_2+\frac{1}{2},m_2,q,p-3>=
  -(J_1+m_1+\frac{7}{2})(J_2+m_2+\frac{7}{2})\n
&&~~~~\times\lt[\frac{1}{2(J_2+1)}
  |J_1,m_1+\frac{1}{2},J_2,m_2+\frac{1}{2},q;p-2>_{\bf I}\rt.\n
&&~~~~\lt.  +\frac{1}{2(J_1+1)} 
  |J_1,m_1+\frac{1}{2},J_2,m_2+\frac{1}{2},q;p-2>_{\bf II}\rt]\n
&&~~~~-\frac{J_2+m_2+\frac{7}{2}}{2(J_1+1)}
  |J_1+1,m_1+\frac{1}{2},J_2,m_2+\frac{1}{2},q;p-2>\n
&&~~~~-\frac{J_1+m_1+\frac{7}{2}}{2(J_2+1)}
  |J_1,m_1+\frac{1}{2},J_2+1,m_2+\frac{1}{2},q;p-2>,\n
&&\G(E_{32})|J_1-\frac{1}{2},m_1,J_2-\frac{1}{2},m_2,q,p-3>\n
&&~~=(q-J_1+J_2)|J_1,m_1-\frac{1}{2},J_2,m_2-\frac{1}{2},q;p-4>,\n
&&\G(E_{32})|J_1+\frac{1}{2},m_1,J_2-\frac{1}{2},m_2,q,p-3>\n
&&~~=(q+J_1+J_2+1)(J_1-m_1-\frac{5}{2})
  |J_1,m_1-\frac{1}{2},J_2,m_2-\frac{1}{2},q;p-4>,\n
&&\G(E_{32})|J_1-\frac{1}{2},m_1,J_2+\frac{1}{2},m_2,q,p-3>\n
&&~~=(q-J_1-J_2-1)(J_2-m_2-\frac{5}{2})
  |J_1,m_1-\frac{1}{2},J_2,m_2-\frac{1}{2},q;p-4>,\n
&&\G(E_{32})|J_1+\frac{1}{2},m_1,J_2+\frac{1}{2},m_2,q,p-3>\n
&&~~=(q+J_1-J_2)(J_1-m_1-\frac{5}{2})(J_2-m_2-\frac{5}{2})
  |J_1,m_1-\frac{1}{2},J_2,m_2-\frac{1}{2},q;p-4>.\n\label{action-level3}
\eeqa

Finally, the actions of the odd simple generators on the level-4
multiplet are
\beqa
&&\G(E_{23})|J_1,m_1,J_2,m_2,q;p-4>=\frac{1}{(2J_1+1)(2J_2+1)}\n
&&~~~~\times\lt[(J_1+m_1+4)(J_2+m_2+4)
  |J_1-\frac{1}{2},m_1+\frac{1}{2},J_2-\frac{1}{2},m_2+\frac{1}{2},q,p-3>\rt.\n
&&~~~~+(J_2+m_2+4)
  |J_1+\frac{1}{2},m_1+\frac{1}{2},J_2-\frac{1}{2},m_2+\frac{1}{2},q,p-3>\n
&&~~~~+(J_1+m_1+4)
  |J_1-\frac{1}{2},m_1+\frac{1}{2},J_2+\frac{1}{2},m_2+\frac{1}{2},q,p-3>\n
&&~~~~\lt.+|J_1+\frac{1}{2},m_1+\frac{1}{2},J_2+\frac{1}{2},
  m_2+\frac{1}{2},q,p-3>\rt],\n
&&\G(E_{32})|J_1,m_1,J_2,m_2,q;p-4>=0.\label{action-level4}
\eeqa

Summarizing, we have obtained 16 independent multiplets,
(\ref{level-04b}) and (\ref{level-1b})--(\ref{level-2b}), of $gl(2)\oplus
gl(2)$ which span finite-dimensional representations of $gl(2|2)$.
For generic $q$, these multiplets span irreducible typical
representations of $gl(2|2)$ of dimension $16(2J_1+1)(2J_2+1)$.
Denote by $\pi_{(J_1,J_2,q,p)}$ and $\sigma_{(J_1,J_2,q,p)}$
the $gl(2|2)$ and $gl(2)\oplus gl(2)$ representations 
with highest weight $(J_1,J_2,q,p)$, respectively.
Then the $gl(2|2)\downarrow gl(2)\oplus gl(2)$ branching rule for generic
$q$ is given by
\beqa
\pi_{(J_1,J_2,q,p)}&=&\sigma_{(J_1,J_2,q,p)}\oplus
\sigma_{(J_1-1/2,J_2-1/2,q,p-1)}\oplus\sigma_{(J_1+1/2,J_2-1/2,q,p-1)}\n
& &\oplus \sigma_{(J_1+1/2,J_2+1/2,q,p-1)}\oplus
\sigma_{(J_1-1/2,J_2+1/2,q,p-1)} \oplus \sigma_{(J_1,J_2-1,q,p-2)}\n
& &\oplus \sigma_{(J_1-1,J_2,q,p-2)} \oplus \sigma_{(J_1+1,J_2,q,p-2)}
\oplus \sigma_{(J_1,J_2+1,q,p-2)} \oplus 2\times \sigma_{(J_1,J_2,q,p-2)}\n
& &\oplus\sigma_{(J_1-1/2,J_2-1/2,q,p-3)}\oplus\sigma_{(J_1+1/2,J_2-1/2,q,p-3)}
\oplus \sigma_{(J_1-1/2,J_2+1/2,q,p-3)}\n
& &\oplus\sigma_{(J_1+1/2,J_2+1/2,q,p-3)} \oplus\sigma_{(J_1,J_2,q,p-4)}.
\eeqa

Some remarks are in order. Firstly,
ireducible representations are obtained as submodules
(not subquotients) of the super-Fock space generated by $\{a_{ij},
a_{ij}^\dg, \a_{ij},\a_{ij}^\dg\}$. This is because the $gl(2|2)$-module
structure of the super-Fock space is the contragredient dual of the
Verma model over $gl(2|2)$. Secondly, as $|J_1,m_1,J_2,m_2,q;p-2>_{\bf I}\equiv
0$ when $J_2=0$ and $|J_1,m_1,J_2,m_2,q;p-2>_{\bf II}\equiv 0$ when
$J_1=0$, thus if $J_1=0$ or $J_2=0$ only one copy of 
$\sigma_{(J_1,J_2,q,p-2)}$ remains in the above branching rule. In
particular, when $J_1=0=J_2$ which corresponds to the 16-dimensional typical
representation of $gl(2|2)$,  $\sigma_{(J_1,J_2,q,p-2)}$
disappears and the branching rule becomes
\beqa
\pi_{(0,0,q,p)}&=&\sigma_{(0,0,q,p)}\oplus\sigma_{(1/2,1/2,q,p-1)}
  \oplus\sigma_{(1,0,q,p-2)}\n
& &\oplus\sigma_{(0,1,q,p-2)}\oplus\sigma_{(1/2,1/2,q,p-3)}
  \oplus\sigma_{(0,0,q,p-4)}
\eeqa
or $\underline{16}=\underline{1}\oplus\underline{4}\oplus\underline{3}
\oplus\underline{3}\oplus\underline{4}\oplus\underline{1}$.

\sect{Atypical Representations of $gl(2|2)$}

We have different types of atypical representations of $gl(2|2)$.
{}From the actions of the odd generators on the $gl(2)\oplus gl(2)$
multiplets, we see that when $q=\pm(J_1-J_2), \pm(J_1+J_2+1)$, the
representations become atypical. The Casimir for such representations 
vanishes, and yet they are not the trivial one-dimensional
representation. 

\subsection{Atypical representation corresponding to $q=J_1-J_2$}

\noindent\underline{Case 1. $q=J_1-J_2,~J_1\neq J_2$}:
\vskip.1in
Let us introduce the following independent combinations:
\beqa
|J_1,m_1,J_2,m_2,q,p-2>_{\bf sym1}&=&
   J_1|J_1,m_1,J_2,m_2,q,p-2>_{\bf I}\n
& &   +J_2|J_1,m_1,J_2,m_2,q,p-2>_{\bf II},\n
|J_1,m_1,J_2,m_2,q,p-2>_{\bf asym1}&=&
   J_1|J_1,m_1,J_2,m_2,q,p-2>_{\bf I}\n
& &   -J_2|J_1,m_1,J_2,m_2,q,p-2>_{\bf II}
\eeqa
for $J_1\neq 0, J_2\neq 0$. When $J_1=0$ or $J_2=0$, we let
$|J_1,m_1,J_2,m_2,q,p-2>_{\bf sym1}\equiv 0$ and
\beq
|J_1,m_1,J_2,m_2,q,p-2>_{\bf asym1}=\lt\{
\begin{array}{ll}
|J_1,m_1,J_2,m_2,q,p-2>_{\bf I}& {\rm if}~ J_1=0,\\
|J_1,m_1,J_2,m_2,q,p-2>_{\bf II}& {\rm if}~ J_2=0.
\end{array}
\rt.
\eeq

It can be shown from the actions of odd generators that when $q=J_1-J_2$, 
\beqa
&&\G(E_{23})|J_1,m_1,J_2,m_2,q,p-2>_{\bf sym1}
  =\frac{1}{(2J_1+1)(2J_2+1)}\times\n
&&~~~~\lt[(J_1-J_2)(J_1+m_1+2)(J_2+m_2+2)|J_1-\frac{1}{2},m_1+\frac{1}{2},
   J_2-\frac{1}{2},m_2+\frac{1}{2},q;p-1>\rt.\n
&&~~~~-J_1(2J_2+1)(J_2+m_2+2)|J_1+\frac{1}{2},m_1+\frac{1}{2},
   J_2-\frac{1}{2},m_2+\frac{1}{2},q;p-1>\n
&&~~~~\lt.-(2J_1+1)J_2(J_1+m_1+2)|J_1-\frac{1}{2},m_1+\frac{1}{2},
   J_2+\frac{1}{2},m_2+\frac{1}{2},q;p-1>\rt]
\eeqa
which does not contain the multiplet 
$|J_1+\frac{1}{2},m_1,J_2+\frac{1}{2},m_2,q;p-1>$ and
\beqa
&&\G(E_{32})|J_1,m_1,J_2,m_2,q,p-2>_{\bf sym1}=\frac{(J_1-J_2)
  (4J_1J_2+2J_1+2J_2+1)}{(2J_1+1)(2J_2+1)}\n
&&~~~~\times
  (J_1-m_1-2)(J_2-m_2-2)|J_1-\frac{1}{2},m_1-\frac{1}{2},J_2-\frac{1}{2},
  m_2-\frac{1}{2},q;p-3>.
\eeqa
Thus when $q=J_1-J_2$, if one starts with the level-0 state
$|J_1,m_1,J_2,m_2,q;p>$
then we find using the actions (\ref{action-level0}-\ref{action-level4})
that the following $gl(2)\oplus gl(2)$ multiplets
\beqa
&&|J_1+\frac{1}{2},m_1,J_2+\frac{1}{2},m_2,q;p-1>,\n
&&|J_1+1,m_1,J_2,m_2,q,p-2>,~~~~|J_1,m_1,J_2+1,m_2,q,p-2>,\n
&&|J_1,m_1,J_2,m_2,q,p-2>_{\bf
  asym1},~~~~|J_1+\frac{1}{2},m_1,J_2-\frac{1}{2},m_2,q;p-3>,\n
&&|J_1-\frac{1}{2},m_1,J_2+\frac{1}{2},m_2,q;p-3>,~~~~
  |J_1+\frac{1}{2},m_1,J_2+\frac{1}{2},m_2,q;p-3>,\n
&&|J_1,m_1,J_2,m_2,q,p-4>
\eeqa
disappear, and only the following multiplets
\beqa
&&|J_1,m_1,J_2,m_2,q,p>,\n
&&|J_1-\frac{1}{2},m_1,J_2-\frac{1}{2},m_2,q;p-1>,~~~~
  |J_1+\frac{1}{2},m_1,J_2-\frac{1}{2},m_2,q;p-1>,\n
&&|J_1-\frac{1}{2},m_1,J_2+\frac{1}{2},m_2,q;p-1>,\n
&&|J_1,m_1,J_2,m_2,q,p-2>_{\bf sym1},~~~~
  |J_1-1,m_1,J_2,m_2,q,p-2>,\n
&&|J_1,m_1,J_2-1,m_2,q,p-2>,~~~~
  |J_1-\frac{1}{2},m_1,J_2-\frac{1}{2},m_2,q;p-3>\label{atypical1}
\eeqa
remain. They form irreducible atypical representations of $gl(2|2)$ of
dimension $8[(2J_1+1)J_2+J_1(2J_2+1)]$. So the $gl(2|2)\downarrow gl(2)\oplus
gl(2)$ branching rule for $q=J_1-J_2$ is given by
\beqa
\pi_{(J_1,J_2,q,p)}&=&\sigma_{(J_1,J_2,q,p)}\oplus\sigma_{(J_1-\frac{1}{2},
  J_2-\frac{1}{2},q,p-1)}\oplus\sigma_{(J_1+\frac{1}{2},J_2+\frac{1}{2},
  q,p-1)}\n
& &\oplus\sigma_{(J_1-\frac{1}{2},J_2+\frac{1}{2},q,p-1)}\oplus
  \sigma_{(J_1,J_2,q,p-2)}\oplus\sigma_{(J_1-1,J_2,q,p-2)}\n
& &\oplus\sigma_{(J_1,J_2-1,q,p-2)}\oplus\sigma_{(J_1-\frac{1}{2},
  J_2-\frac{1}{2},q,p-3)}.
\eeqa
It should be understood here that $\sigma_{(J_1,J_2,q,p-2)}$ disappears
when $J_1=0$ or $J_2=0$.

\vskip.2in
\noindent\underline{Case 2. $q=J_1-J_2,~ J_1=J_2$ so that $q=0$}:
\vskip.1in

In this case, we define the independent combinations:
\beqa
|J_1,m_1,J_2,m_2,q,p-2>_{\bf sym1'}&=&
   |J_1,m_1,J_2,m_2,q,p-2>_{\bf I}\n
& &   +|J_1,m_1,J_2,m_2,q,p-2>_{\bf II},\n
|J_1,m_1,J_2,m_2,q,p-2>_{\bf asym1'}&=&
   |J_1,m_1,J_2,m_2,q,p-2>_{\bf I}\n
& &   -|J_1,m_1,J_2,m_2,q,p-2>_{\bf II}.
\eeqa
Both $|J_1,m_1,J_2,m_2,q,p-2>_{\bf sym1'}$ and
$|J_1,m_1,J_2,m_2,q,p-2>_{\bf asym1'}$ vanish if $J_1=0=J_2$. Then it is
easily shown that $\G(E_{23})
|J_1,m_1,J_2,m_2,q,p-2>_{\bf sym1'}$ does not contain
$|J_1-\frac{1}{2},m_1,J_2-\frac{1}{2},m_2,q;p-1>$ and
$|J_1+\frac{1}{2},m_1,J_2+\frac{1}{2},m_2,q;p-1>$, and
$\G(E_{32})|J_1,m_1,J_2,m_2,q,p-2>_{\bf sym1'}=0$. Thus only the
following multiplets
\beqa
&&|J_1,m_1,J_2,m_2,q,p>,\n
&& |J_1+\frac{1}{2},m_1,J_2-\frac{1}{2},m_2,q;p-1>,~~~~
  |J_1-\frac{1}{2},m_1,J_2+\frac{1}{2},m_2,q;p-1>,\n
&&|J_1,m_1,J_2,m_2,q,p-2>_{\bf sym1'}\label{j1=j2}
\eeqa
survive, and they give irreducible atypical representations of dimension
$4[(2J_1+1)(2J_2+1)-1/2]$ if $J_1=J_2\neq 0$ and the trivial
one-deimsional representaion if $J_1=0=J_2$ (for which the last
three multiplets in (\ref{j1=j2}) disappear).

\vskip.2in
\noindent\underline {Case 3. Lowest weight (indecomposable) Kac modules}:
\vskip.1in
Other types of atypical representations when $q=J_1-J_2$  are not
irreducible. One such type of representations are obtained by starting with
the level-4 state $|J_1,m_1,J_2,m_2,q;p-4>$.
These representations contain all 16 multiplets and a non-separable
invariant subspace provided by the multiplets (\ref{atypical1})
[or (\ref{j1=j2}) when $J_1=J_2$]. These representations are
not fully reducible (i.e. indecomposable) and have dimension
$16(2J_1+1)(2J_2+1)$.

\subsection{Atypical representations corresponding to $q=-J_1+J_2$}

The case where $J_1=J_2$ so that $q=0$ is the same as Case 2 of the
last subsection. So in this subsection we only consider the $J_1\neq J_2$
case.

\vskip.2in
\noindent\underline{1. Irreducible representations}:
\vskip.1in
Let us introduce the following independent combinations:
\beqa
|J_1,m_1,J_2,m_2,q,p-2>_{\bf sym2}&=&
   (J_1+1)|J_1,m_1,J_2,m_2,q,p-2>_{\bf I}\n
& &   +(J_2+1)|J_1,m_1,J_2,m_2,q,p-2>_{\bf II},\n
|J_1,m_1,J_2,m_2,q,p-2>_{\bf asym2}&=&
   (J_1+1)|J_1,m_1,J_2,m_2,q,p-2>_{\bf I}\n
& &   -(J_2+1)|J_1,m_1,J_2,m_2,q,p-2>_{\bf II}
\eeqa
for $J_1\neq 0, J_2\neq 0$, and let 
\beq
|J_1,m_1,J_2,m_2,q,p-2>_{\bf sym2}=\lt\{
\begin{array}{ll}
|J_1,m_1,J_2,m_2,q,p-2>_{\bf I}& {\rm if}~ J_1=0,\\
|J_1,m_1,J_2,m_2,q,p-2>_{\bf II}& {\rm if}~ J_2=0.
\end{array}
\rt.
\eeq
and $|J_1,m_1,J_2,m_2,q,p-2>_{\bf asym2}=0$  if 
$J_1=0$ or $J_2=0$. 

Similar to the $q=J_1-J_2$ case, we may show that when $q=-J_1+J_2$, 
\beqa
&&\G(E_{23})|J_1,m_1,J_2,m_2,q,p-2>_{\bf sym2}
  =\frac{1}{(2J_1+1)(2J_2+1)}\times\n
&&~~~~\lt[-(2J_1+1)(J_2+1)(J_2+m_2+2)|J_1+\frac{1}{2},m_1+\frac{1}{2},
   J_2-\frac{1}{2},m_2+\frac{1}{2},q;p-1>\rt.\n
&&~~~~+(J_1-J_2)|J_1+\frac{1}{2},m_1+\frac{1}{2},
   J_2+\frac{1}{2},m_2+\frac{1}{2},q;p-1>\n
&&~~~~\lt.-(J_1+1)(2J_2+1)(J_1+m_1+2)|J_1-\frac{1}{2},m_1+\frac{1}{2},
   J_2+\frac{1}{2},m_2+\frac{1}{2},q;p-1>\rt]\n
\eeqa
which is independent of 
$|J_1-\frac{1}{2},m_1,J_2-\frac{1}{2},m_2,q;p-1>$ and
\beqa
&&\G(E_{32})|J_1,m_1,J_2,m_2,q,p-2>_{\bf sym2}=\frac{(J_1-J_2)
  (4J_1J_2+2J_1+2J_2+1)}{(2J_1+1)(2J_2+1)}\n
&&~~~~\times
  (J_1-m_1-2)(J_2-m_2-2)|J_1+\frac{1}{2},m_1+\frac{1}{2},J_2+\frac{1}{2},
  m_2+\frac{1}{2},q;p-3>.
\eeqa
Thus when $q=-J_1+J_2$, if one starts with the level-0
state then by the actions (\ref{action-level0}-\ref{action-level4})
one finds that the following $gl(2)\oplus gl(2)$ multiplets
\beqa
&&|J_1-\frac{1}{2},m_1,J_2-\frac{1}{2},m_2,q;p-1>,\n
&&|J_1,m_1,J_2-1,m_2,q,p-2>,~~~~|J_1-1,m_1,J_2,m_2,q,p-2>,\n
&&|J_1,m_1,J_2,m_2,q,p-2>_{\bf
  asym2},~~~~|J_1-\frac{1}{2},m_1,J_2-\frac{1}{2},m_2,q;p-3>,\n
&&|J_1+\frac{1}{2},m_1,J_2-\frac{1}{2},m_2,q;p-3>,~~~~
  |J_1-\frac{1}{2},m_1,J_2+\frac{1}{2},m_2,q;p-3>,\n
&&|J_1,m_1,J_2,m_2,q,p-4>
\eeqa
drop out of the basis, and only the following multiplets
\beqa
&&|J_1,m_1,J_2,m_2,q,p>,\n
&&|J_1+\frac{1}{2},m_1,J_2+\frac{1}{2},m_2,q;p-1>,~~~~
  |J_1+\frac{1}{2},m_1,J_2-\frac{1}{2},m_2,q;p-1>,\n
&&|J_1-\frac{1}{2},m_1,J_2+\frac{1}{2},m_2,q;p-1>,\n
&&|J_1,m_1,J_2,m_2,q,p-2>_{\bf sym2},~~~~
  |J_1+1,m_1,J_2,m_2,q,p-2>,\n
&&|J_1,m_1,J_2+1,m_2,q,p-2>,~~~~
  |J_1+\frac{1}{2},m_1,J_2+\frac{1}{2},m_2,q;p-3>\label{atypical2}
\eeqa
survive. They form irreducible
atypical representations of $gl(2|2)$ of dimension
$8[(J_1+1)(2J_2+1)+(2J_1+1)(J_2+1)]$. The branching rule in this case 
(i.e. $q=-J_1+J_2$) becomes
\beqa
\pi_{(J_1,J_2,q,p)}&=&\sigma_{(J_1,J_2,q,p)}\oplus\sigma_{(J_1+\frac{1}{2},
  J_2+\frac{1}{2},q,p-1)}\oplus\sigma_{(J_1+\frac{1}{2},J_2-\frac{1}{2},
  q,p-1)}\n
& &\oplus\sigma_{(J_1-\frac{1}{2},J_2+\frac{1}{2},q,p-1)}\oplus
  \sigma_{(J_1,J_2,q,p-2)}\oplus\sigma_{(J_1+1,J_2,q,p-2)}\n
& &\oplus\sigma_{(J_1,J_2+1,q,p-2)}\oplus\sigma_{(J_1+\frac{1}{2},
  J_2+\frac{1}{2},q,p-3)}.
\eeqa

\vskip.1in
\noindent\underline{2. Lowest weight (indecomposable) Kac modules}:
\vskip.1in

If one starts with the level-4 state, then one gets
atypical representations which are not irreducible.
In such representations, all 16 multiplets appear but there exists
a non-separable invariant superspace
generated by multiplets (\ref{atypical2}). These representations are
indecomposable and have dimension $16(2J_1+1)(2J_2+1)$.

\subsection{Atypical representations corresponding to $q=J_1+J_2+1$}

\vskip.1in
\noindent\underline{1. Irreducible representations}:
\vskip.1in
Let us introduce the following independent combinations for $J_1\neq 0,
J_2\neq 0$,
\beqa
|J_1,m_1,J_2,m_2,q,p-2>_{\bf sym3}&=&
   J_1\,|J_1,m_1,J_2,m_2,q,p-2>_{\bf I}\n
& &   +(J_2+1)|J_1,m_1,J_2,m_2,q,p-2>_{\bf II},\n
|J_1,m_1,J_2,m_2,q,p-2>_{\bf asym3}&=&
   J_1\,|J_1,m_1,J_2,m_2,q,p-2>_{\bf I}\n
& &   -(J_2+1)|J_1,m_1,J_2,m_2,q,p-2>_{\bf II}.
\eeqa
We let
\beqa
|J_1,m_1,J_2,m_2,q,p-2>_{\bf sym3}&=&\lt\{
\begin{array}{ll}
|J_1,m_1,J_2,m_2,q,p-2>_{\bf I}& {\rm if}~ J_1=0,\\
0 & {\rm if}~ J_2=0,
\end{array}
\rt.\n
|J_1,m_1,J_2,m_2,q,p-2>_{\bf asym3}&=&\lt\{
\begin{array}{ll}
0 & {\rm if}~ J_1=0,\\
|J_1,m_1,J_2,m_2,q,p-2>_{\bf II}& {\rm if}~ J_2=0.
\end{array}
\rt.
\eeqa

It can be seen from the actions of odd generators that when $q=J_1+J_2+1$, 
\beqa
&&\G(E_{23})|J_1,m_1,J_2,m_2,q,p-2>_{\bf asym3}
  =\lt.\frac{1}{(2J_1+1)(2J_2+1)}\rt[(2J_1+1)(J_2+1)\n
&&~~~~\times(J_2+m_2+2)(J_2+m_2+2)|J_1-\frac{1}{2},m_1+\frac{1}{2},
   J_2-\frac{1}{2},m_2+\frac{1}{2},q;p-1>\n
&&~~~~-J_1(2J_2+1)|J_1+\frac{1}{2},m_1+\frac{1}{2},
   J_2+\frac{1}{2},m_2+\frac{1}{2},q;p-1>\n
&&~~~~\lt.+(J_1+J_2+1)(J_1+m_1+2)|J_1-\frac{1}{2},m_1+\frac{1}{2},
   J_2+\frac{1}{2},m_2+\frac{1}{2},q;p-1>\rt]
\eeqa
which does not contain the multiplet 
$|J_1+\frac{1}{2},m_1,J_2-\frac{1}{2},m_2,q;p-1>$ and
\beqa
&&\G(E_{32})|J_1,m_1,J_2,m_2,q,p-2>_{\bf asym3}=\frac{(J_1+J_2+1)
  (4J_1J_2+J_1+J_2)+(J_1+J_2)^2}{(2J_1+1)(2J_2+1)}\n
&&~~~~\times
  (J_1-m_1-2)|J_1-\frac{1}{2},m_1-\frac{1}{2},J_2+\frac{1}{2},
  m_2-\frac{1}{2},q;p-3>.
\eeqa
Then similar to previous cases, when $q=J_1+J_2+1$,
the following $gl(2)\oplus gl(2)$ multiplets
\beqa
&&|J_1+\frac{1}{2},m_1,J_2-\frac{1}{2},m_2,q;p-1>,\n
&&|J_1,m_1,J_2-1,m_2,q,p-2>,~~~~|J_1+1,m_1,J_2,m_2,q,p-2>,\n
&&|J_1,m_1,J_2,m_2,q,p-2>_{\bf
  sym3},~~~~|J_1-\frac{1}{2},m_1,J_2-\frac{1}{2},m_2,q;p-3>,\n
&&|J_1+\frac{1}{2},m_1,J_2-\frac{1}{2},m_2,q;p-3>,~~~~
  |J_1+\frac{1}{2},m_1,J_2+\frac{1}{2},m_2,q;p-3>,\n
&&|J_1,m_1,J_2,m_2,q,p-4>
\eeqa
disappear, and only the following multiplets
\beqa
&&|J_1,m_1,J_2,m_2,q,p>,\n
&&|J_1-\frac{1}{2},m_1,J_2-\frac{1}{2},m_2,q;p-1>,~~~~
  |J_1+\frac{1}{2},m_1,J_2+\frac{1}{2},m_2,q;p-1>,\n
&&|J_1-\frac{1}{2},m_1,J_2+\frac{1}{2},m_2,q;p-1>,\n
&&|J_1,m_1,J_2,m_2,q,p-2>_{\bf asym3},~~~~
  |J_1-1,m_1,J_2,m_2,q,p-2>,\n
&&|J_1,m_1,J_2+1,m_2,q,p-2>,~~~~
  |J_1-\frac{1}{2},m_1,J_2+\frac{1}{2},m_2,q;p-3>\label{atypical3}
\eeqa
remain. They constitute irreducible
atypical representations of $gl(2|2)$ of dimension
$8[(2J_1+1)(J_2+1)+J_1(2J_2+1)]$. The branching rule in this case 
(i.e. $q=J_1+J_2+1$) reads
\beqa
\pi_{(J_1,J_2,q,p)}&=&\sigma_{(J_1,J_2,q,p)}\oplus\sigma_{(J_1-\frac{1}{2},
  J_2-\frac{1}{2},q,p-1)}\oplus\sigma_{(J_1+\frac{1}{2},J_2+\frac{1}{2},
  q,p-1)}\n
& &\oplus\sigma_{(J_1-\frac{1}{2},J_2+\frac{1}{2},q,p-1)}\oplus
  \sigma_{(J_1,J_2,q,p-2)}\oplus\sigma_{(J_1-1,J_2,q,p-2)}\n
& &\oplus\sigma_{(J_1,J_2+1,q,p-2)}\oplus\sigma_{(J_1-\frac{1}{2},
  J_2+\frac{1}{2},q,p-3)}.
\eeqa
Here one should keep in mind that
$\sigma_{(J_1,J_2,q,p-2)}$ disappears if $J_1=0$.

\vskip.2in
\noindent\underline{2. Lowest weight (indecomposable) Kac representations}:
\vskip.1in
Similar to the previous cases, if one retains all 16 multiplets, then
one gets lowest weight (indecomposable) Kac representations of 
$16(2J_1+1)(2J_2+1)$ which contain an invariant but non-separable subspace
provided by multiplets (\ref{atypical3}).

\subsection{Atypical representations corresponding to $q=-J_1-J_2-1$}

\vskip.1in
\noindent\underline{1. Irreducible representations}:
\vskip.1in
In this case, we introduce the following independent combinations for 
$J_1\neq 0, J_2\neq 0$,
\beqa
|J_1,m_1,J_2,m_2,q,p-2>_{\bf sym4}&=&
   (J_1+1)|J_1,m_1,J_2,m_2,q,p-2>_{\bf I}\n
& &   +J_2\,|J_1,m_1,J_2,m_2,q,p-2>_{\bf II},\n
|J_1,m_1,J_2,m_2,q,p-2>_{\bf asym4}&=&
   (J_1+1)|J_1,m_1,J_2,m_2,q,p-2>_{\bf I}\n
& &   -J_2\,|J_1,m_1,J_2,m_2,q,p-2>_{\bf II}
\eeqa
and let
\beqa
|J_1,m_1,J_2,m_2,q,p-2>_{\bf sym4}&=&\lt\{
\begin{array}{ll}
0& {\rm if}~ J_1=0,\\
|J_1,m_1,J_2,m_2,q,p-2>_{\bf II} & {\rm if}~ J_2=0,
\end{array}
\rt.\n
|J_1,m_1,J_2,m_2,q,p-2>_{\bf asym4}&=&\lt\{
\begin{array}{ll}
|J_1,m_1,J_2,m_2,q,p-2>_{\bf I} & {\rm if}~ J_1=0,\\
0& {\rm if}~ J_2=0.
\end{array}
\rt.
\eeqa

It can be seen from the actions of odd generators that when $q=-J_1-J_2-1$, 
\beqa
&&\G(E_{23})|J_1,m_1,J_2,m_2,q,p-2>_{\bf asym4}
  =\lt.\frac{1}{(2J_1+1)(2J_2+1)}\rt[(J_1+1)(2J_2+1)\n
&&~~~~\times(J_2+m_2+2)(J_2+m_2+2)|J_1-\frac{1}{2},m_1+\frac{1}{2},
   J_2-\frac{1}{2},m_2+\frac{1}{2},q;p-1>\n
&&~~~~-(J_1+J_2+1)(J_2+m_2+2)|J_1+\frac{1}{2},m_1+\frac{1}{2},
   J_2-\frac{1}{2},m_2+\frac{1}{2},q;p-1>\n
&&~~~~\lt.-(2J_1+1)J_2|J_1+\frac{1}{2},m_1+\frac{1}{2},
   J_2+\frac{1}{2},m_2+\frac{1}{2},q;p-1>\rt]
\eeqa
which has no dependence on the multiplet 
$|J_1-\frac{1}{2},m_1,J_2+\frac{1}{2},m_2,q;p-1>$ and
\beqa
&&\G(E_{32})|J_1,m_1,J_2,m_2,q,p-2>_{\bf asym4}=-\frac{(J_1+J_2+1)
  (4J_1J_2+J_1+J_2)+(J_1+J_2)^2}{(2J_1+1)(2J_2+1)}\n
&&~~~~\times
  (J_2-m_2-2)|J_1+\frac{1}{2},m_1-\frac{1}{2},J_2-\frac{1}{2},
  m_2-\frac{1}{2},q;p-3>.
\eeqa
Thus when $q=-J_1-J_2-1$, the following $gl(2)\oplus gl(2)$ multiplets
\beqa
&&|J_1-\frac{1}{2},m_1,J_2+\frac{1}{2},m_2,q;p-1>,\n
&&|J_1-1,m_1,J_2,m_2,q,p-2>,~~~~|J_1,m_1,J_2+1,m_2,q,p-2>,\n
&&|J_1,m_1,J_2,m_2,q,p-2>_{\bf
  sym4},~~~~|J_1-\frac{1}{2},m_1,J_2-\frac{1}{2},m_2,q;p-3>,\n
&&|J_1-\frac{1}{2},m_1,J_2+\frac{1}{2},m_2,q;p-3>,~~~~
  |J_1+\frac{1}{2},m_1,J_2+\frac{1}{2},m_2,q;p-3>,\n
&&|J_1,m_1,J_2,m_2,q,p-4>
\eeqa
drop out, and only the following multiplets
\beqa
&&|J_1,m_1,J_2,m_2,q,p>,\n
&&|J_1-\frac{1}{2},m_1,J_2-\frac{1}{2},m_2,q;p-1>,~~~~
  |J_1+\frac{1}{2},m_1,J_2+\frac{1}{2},m_2,q;p-1>,\n
&&|J_1+\frac{1}{2},m_1,J_2-\frac{1}{2},m_2,q;p-1>,\n
&&|J_1,m_1,J_2,m_2,q,p-2>_{\bf asym4},~~~~
  |J_1,m_1,J_2-1,m_2,q,p-2>,\n
&&|J_1+1,m_1,J_2,m_2,q,p-2>,~~~~
  |J_1+\frac{1}{2},m_1,J_2-\frac{1}{2},m_2,q;p-3>\label{atypical4}
\eeqa
remain. They give irreducible atypical representations of $gl(2|2)$ of dimension
$8[(J_1+1)(2J_2+1)+(2J_1+1)J_2]$. In this case the branching rule becomes
\beqa
\pi_{(J_1,J_2,q,p)}&=&\sigma_{(J_1,J_2,q,p)}\oplus\sigma_{(J_1-\frac{1}{2},
  J_2-\frac{1}{2},q,p-1)}\oplus\sigma_{(J_1+\frac{1}{2},J_2+\frac{1}{2},
  q,p-1)}\n
& &\oplus\sigma_{(J_1+\frac{1}{2},J_2-\frac{1}{2},q,p-1)}\oplus
  \sigma_{(J_1,J_2,q,p-2)}\oplus\sigma_{(J_1,J_2-1,q,p-2)}\n
& &\oplus\sigma_{(J_1+1,J_2,q,p-2)}\oplus\sigma_{(J_1+\frac{1}{2},
  J_2-\frac{1}{2},q,p-3)}.
\eeqa
Here it should be understood that $\sigma_{(J_1,J_2,q,p-2)}$ is not in
the branching rule if $J_2=0$.

\vskip.2in
\noindent\underline{2. Lowest weight (indecomposable) Kac representations}:
\vskip.1in
As before, other types of atypical representations are not irreducible.
These representations contain all 16 multiplets which contain a
non-separable invariant subspace generated by multiplets
(\ref{atypical4}). They are lowest weight (indecomposable) Kac
representations of dimension $16(2J_1+1)(2J_2+1)$.

\section{Conclusions and Discussions}

In this article we have applied the super coherent
state method to the construction of the free boson-fermion realization and
representations of the non-semisimple superalgebra $gl(2|2)$  in the
standard basis. The representations are constructed out of the 
$gl(2)\oplus gl(2)$ particle states in the super-Fock space. 

As mentioned in the introduction, 
superalgebras and their corresponding non-unitary CFTs emerge in the
supersymmetric treatment to disordered systems and the integer quantum
Hall plateaus. In such a treatment, primary fields play an important role
in the computation of critical properties of the disordered systems.
The results obtained in this paper now
make possible the construction of all primary fields of the $gl(2|2)$
non-unitary CFT in terms of free fields \cite{Din03b}. This is under
investigation and results will be presented elsewhere.

\vskip.3in

\no {\bf Acknowledgments:}
Our interest in the coherent state construction was ignited by Max
Lohe's talk \cite{Loh02}. We thank Max Lohe for making the talk material
available to us.  The financial support from the Australian Research 
Council is gratefully acknowledged.

\bebb{99}

\bbit{Roz92} 
L. Rozanski and H. Saleur, \npb {376} {461} {1992}.

\bbit{Isi94}
J. M. Isidro and A. V. Ramallo, \npb {414} {715} {1994}. 

\bbit{Flohr}
M. Flohr, Int. J. Mod. Phys. {\bf A28}, 4497 (2003).

\bbit{Efe83}
K. Efetov, Adv. Phys. {\bf 32}, 53 (1983). 

\bbit{Ber95}
D. Bernard, preprint hep-th/9509137.

\bbit{Mud96}
C. Mudry, C. Chamon and X.-G. Wen, \npb {466} {383} {1996}.

\bbit{Maa97}
Z. Maassarani and D. Serban, \npb {489} {603} {1997}.

\bbit{Zir99}
M.R. Zirnbauer, preprint hep-th/9905054.

\bbit{Bas00}
Z.S. Bassi and A. LeClair, \npb {578} {577} {2000}.

\bbit{Gur00}
S. Guruswamy, A. LeClair and A.W.W. Ludwig, \npb {583} {475} {2000}.

%\bbit{Lud00}
%A.W.W. Ludwig, preprint cond-mat/0012189.

\bbit{Bha01}
M. J. Bhaseen,  J.-S. Caux ,  I. I. Kogan  and  A. M. Tsvelik, 
\npb {618} {465} {2001}.

%\bbit{Bal89}
%A.B. Balantekin, H.A. Schmitt and P. Halse, \jmp {30} {274} {1989}.

\bbit{Sch77} 
M. Scheunert, W. Nahm and V. Rittenberg, \jmp {18} {155} {1977};
\jmp {18} {146} {1977}.

\bbit{Mar80} 
M. Marcu, \jmp {21} {1277} {1980}; \jmp {21} {1284} {1980}.

%\bbit{Gou89}
%M.D. Gould, J. Phys. {\bf A22}, 1209 (1989).

\bbit{Zha03a}
Y.Z. Zhang, Phys. Lett. {\bf A327}, 442 (2004).

\bbit{Zha03b}
Y.Z. Zhang, hep-th/0405066, to appear in ``Progress in Field Theory Research",
Nova Science Publishers Inc., New York, 2004.

\bbit{Din03a}
X.M. Ding, M. D. Gould, C. J. Mewton and Y. Z. Zhang, 
 J. Phys. {\bf A36}, 7649 (2003).

\bbit{Bow96}
P. Bowcock, R-L.K. Koktava and A. Taormina, \plb {388} {303} {1996}.

\bbit{Ras98}
J. Rasmussen, \npb {510} {688} {1998}.

\bbit{Kam89}
A.H. Kamupingene, N.A. Ky and T.D. Palev, \jmp {30} {553}{1989}.

\bbit{Pal90}
T.D. Palev and N.I. Stoilova, \jmp {31} {953} {1990}.

\bbit{Loh02} M. Lohe, ``Vector coherent states and quantum
affine algebras", talk given at the 3rd University of Queensland
Mathematical Physics Workshop, Oct 2-4, 2002, Coolangatta, Australia.

\bbit{Bal88}
A.B. Balantekin, H.A. Schmitt and B.R. Barrett, \jmp {29} {1634} {1988}.

\bbit{Din03b}
X.M. Ding, M. D. Gould and Y. Z. Zhang,
Phys. Lett. {\bf A318}, 354 (2003).

\eebb

\end{document}